\documentclass[11pt]{article}
\usepackage[margin=1in]{geometry}
\usepackage{url}
\usepackage{booktabs}
\usepackage{amsfonts}
\usepackage{xcolor}
\usepackage{graphicx}
\usepackage{multirow}
\usepackage{color}
\usepackage{amsmath}\allowdisplaybreaks
\usepackage{amsfonts,bm}

\def\1{\bf{1}}

\newcommand{\Norm}[1]{\left\| #1 \right\|}

\def\inner#1#2{\langle #1, #2 \rangle}

\def\fH{{\mathcal{H}}}

\def\VI{{\rm{VI}}}

\DeclareMathOperator*{\argmax}{arg\,max}
\DeclareMathOperator*{\argmin}{arg\,min}

\usepackage{etoolbox}
\usepackage{amsthm}
\theoremstyle{plain}

\makeatletter
\def\Ddots{\mathinner{\mkern1mu\raise\p@
\vbox{\kern7\p@\hbox{.}}\mkern2mu
\raise4\p@\hbox{.}\mkern2mu\raise7\p@\hbox{.}\mkern1mu}}
\makeatother

\makeatletter
\newcommand*{\rom}[1]{\expandafter\@slowromancap\romannumeral #1@}
\makeatother

\newtheorem{theorem}{Theorem}[section]

\newtheorem{lemma}[theorem]{Lemma}

\theoremstyle{definition}
\newtheorem{definition}[theorem]{Definition}
\newtheorem{assumption}[theorem]{Assumption}
\theoremstyle{remark}
\newtheorem{remark}[theorem]{Remark}

\def\oo{{\bf o}}

\def\u{{\bf u}}

\def\w{{\bf w}}

\def\x{{\bf x}}

\def\y{{\bf y}}

\def\z{{\bf z}}
\def\0{{\bf 0}}
\def\1{{\bf 1}}

\def\OM{{\mathcal O}}

\def\ZM{{\mathcal Z}}

\def\RB{{\mathbb R}}

\def\TT{{\rm{T}}}

\usepackage{algorithm}
\usepackage{algorithmic}
\usepackage{enumitem}
\usepackage{hhline}
\def\defeq{\overset{\rm def}{=}}

\def\Holder{{H\oo lder}}
\def \oo{ö}

\usepackage{MnSymbol}
\usepackage{hyperref}
\def\u{{\bf u}}
\AtBeginEnvironment{algorithmic}{\normalsize}

\title{Universal Tensor Methods for Monotone Variational Inequalities}
\author{
\begin{tabular}{c}
Chengchang Liu\\
Department of Artificial Intelligence, Westlake University\\
\texttt{liuchengchang@westlake.edu.cn}\\[0.5em]
John C.S.~Lui\\
Department of Computer Science and Engineering\\
The Chinese University of Hong Kong\\
\texttt{cslui@cse.cuhk.edu.hk}\\[0.5em]
Luo Luo\\
School of Data Science, Fudan University\\
\texttt{luoluo@fudan.edu.cn}
\end{tabular}
}
\date{}

\begin{document}
\maketitle

\begin{abstract}
We study monotone variational inequalities whose operators have \Holder~continuous higher-order derivatives.
For a fixed order $p\geq 2$, we assume that the $(p-1)$-th derivative of the monotone operator is \Holder~continuous with parameter $\nu\in[0,1]$ on a bounded closed convex set.
We develop regularized tensor extragradient methods that combine a high-order Taylor approximation of the operator with an extragradient correction step.
When the \Holder~parameter $\nu$ is known, our regularized tensor extragradient method finds an $\epsilon$-weak solution using $\OM(\epsilon^{-2/(p+\nu)})$ tensor-oracle calls.
When $\nu$ is unknown, we propose a universal tensor extragradient method whose tensor-oracle complexity is $\OM(\epsilon^{-2p/((p+1)(p-1+\nu))})$.
\end{abstract}



\section{Introduction}

We consider the monotone variational inequality
\begin{align}
\label{eq:VI}
\text{find } \z^*\in\ZM \quad \text{such that}\quad
\inner{F(\z^*)}{\z-\z^*}\geq 0,\qquad \text{for all }\z\in\ZM,
\end{align}
where $\ZM\subseteq \RB^d$ is closed and convex, and $F:\ZM\to\RB^d$ is monotone.
We denote the solution set of problem \eqref{eq:VI} as
\begin{align*}
    \VI_\ZM(F)
    \defeq
    \left\{
    \z^*\in\ZM:
    \inner{F(\z^*)}{\z-\z^*}\geq0,\quad \text{for all }\z\in\ZM
    \right\}.
\end{align*}
Variational inequalities are widely studied in economics and game theory~\cite{bacsar1998dynamic,von2007theory}, and they also appear in machine learning applications such as robust optimization~\cite{ben2009robust} and generative adversarial learning~\cite{goodfellow2014generative}.
They form a broad optimization model that includes constrained convex minimization, convex-concave saddle-point problems, and fixed-point problems.
For background on finite-dimensional variational inequalities, we refer the reader to Facchinei and Pang~\cite{facchinei2003finite}; for a recent survey of smooth monotone VI models, see Beznosikov et al.~\cite{beznosikov2022smooth}.

First-order methods for monotone variational inequalities have been studied extensively.
When the operator $F$ is Lipschitz continuous, classical extragradient methods~\cite{korpelevich1976extragradien,tseng1995linear,nemirovski2004prox,popov1980modification}, hybrid proximal extragradient methods~\cite{solodov1999hybrid}, optimistic gradient methods~\cite{DaskalakisISZ18,mokhtari2020convergence,popov1980modification}, and dual extrapolation~\cite{nesterov2006solving} find an $\epsilon$-weak solution in $\OM(\epsilon^{-1})$ iterations.
Second-order methods can improve this dependence when the Jacobian of $F$ is Lipschitz continuous.
Regularized Newton-type methods for monotone VIs were studied by Nesterov and Scrimali~\cite{nesterov2006solving} and Huang et al.~\cite{huang2020cubic}, who established global convergence guarantees and local fast convergence properties.
Monteiro and Svaiter~\cite{monteiro2012iteration} proposed the Newton proximal extragradient method, which achieves the rate $\OM(\epsilon^{-2/3}\log\log(\epsilon^{-1}))$ under Lipschitz Jacobian assumptions.
The extra $\log\log(\epsilon^{-1})$ factor comes from the implicit search used to solve the regularized Newton subproblem.
Subsequent work removed this logarithmic factor: search-free proximal-extragradient schemes~\cite{alves2023search}, second-order extensions of optimistic gradient methods~\cite{jiang2022generalized}, and second-order dual-extrapolation methods~\cite{lin2022perseus,lin2022explicit} achieve the rate $\OM(\epsilon^{-2/3})$ for finding an $\epsilon$-weak solution.

For a general order $p\geq2$, $p$-th order methods for monotone VIs, which use derivatives up to order $p-1$, have been developed~\cite{bullins2022higher,huang2022approximation,adil2022line,nesterov2023high,lin2022perseus}.
Among them, the tensor extragradient methods of Adil et al.~\cite{adil2022line} and Huang and Zhang~\cite{huang2022approximation} are based on iterations of the following form:
\begin{align}
\label{eq:prior_tensor_step}
\begin{cases}
\displaystyle
\z_{k+1/2}
=
\VI_\ZM\left(
\sum_{i=0}^{p-1}\frac{1}{i!}\nabla^iF(\z_k)[\z-\z_k]^i
    +M\|\z-\z_k\|^{p-1}(\z-\z_k)
\right),
\\[0.35cm]
\displaystyle
\z_{k+1}
=
\argmin_{\z\in\ZM}
\left\{
\inner{F(\z_{k+1/2})}{\z\!-\z_{k+1/2}}
+\frac{M\|\z_{k+1/2}-\z_k\|^{p-1}}{2}\|\z-\z_{k+1/2}\|^2
\right\}.
\end{cases}
\end{align}
Under the Lipschitz continuity of $\nabla^{p-1}F$, these methods obtain the rate $\OM(\epsilon^{-2/(p+1)})$ for finding an $\epsilon$-weak solution.

All existing tensor methods for monotone VIs mentioned above are based on the Lipschitz continuity of the highest derivative, while this is not the only possibility for measuring smoothness.
A natural generalization of the Lipschitz continuity assumption is \Holder~continuity.
In convex optimization, Grapiglia and Nesterov~\cite{grapiglia2017regularized,grapiglia2019accelerated,grapiglia2020tensor} proposed regularized Newton-type and tensor algorithms and established the corresponding convergence analysis under the assumption of \Holder~continuous higher-order derivatives.
However, these results cannot be used directly to solve general variational inequality problems because they depend on the convexity of the objective function and on objective-value comparisons.

In this paper, we propose regularized tensor extragradient methods for solving monotone variational inequalities whose $(p-1)$-th derivative is \Holder~continuous with parameter $\nu\in[0,1]$.
We design the algorithms by integrating adaptive search of the regularization parameter with tensor extragradient steps.
The adaptive algorithms for convex minimization typically rely on tracking the objective function value~\cite{grapiglia2017regularized,grapiglia2019accelerated,doikov2021minimizing,cartis2011adaptive,cartis2011adaptive2,nesterov2008accelerating}.
However, for variational inequalities, the weak gap is generally not available during the iterations.
To address this issue, we establish computable criteria based on the operator-model residual to search for the regularization parameter.
The main contributions are as follows.
\begin{itemize}
    \item When the \Holder~parameter $\nu$ is known, we propose the $\nu$-regularized tensor extragradient method ($\nu$-RTEG) and its adaptive variant ($\nu$-ATEG). Both methods find an $\epsilon$-weak solution using $\OM(\epsilon^{-2/(p+\nu)})$ tensor-oracle calls; the adaptive variant does not require the \Holder~constant $\fH_{p,\nu}$.
    \item When neither $\nu$ nor the \Holder~constant is known, we propose a universal tensor extragradient method (UTEG). Its tensor-oracle complexity for finding an $\epsilon$-weak solution is $\OM(\epsilon^{-2p/((p+1)(p-1+\nu))})$.
\end{itemize}
We compare our results with the existing methods in Table~\ref{table:tensor}.
\begin{table}[t]
\caption{Comparison of tensor methods for finding an $\epsilon$-weak solution of a monotone variational inequality with a \Holder~continuous $(p-1)$-th derivative, i.e., $\|\nabla^{p-1}F(\x)-\nabla^{p-1}F(\y)\|\leq \fH_{p,\nu}\|\x-\y\|^\nu$.}
\label{table:tensor}
\centering
\scriptsize
\renewcommand{\arraystretch}{1.15}
\setlength{\tabcolsep}{2.4pt}
\begin{tabular}{@{}lcccp{0.18\textwidth}@{}}
\toprule
Method & Complexity & \Holder~parameter & Required parameter & Reference \\
\midrule
Tensor VI & $\OM(\epsilon^{-\frac{2}{p+1}})$ & $\nu=1$ & $\fH_{p,1}$ & \cite{adil2022line,bullins2022higher,huang2022approximation,lin2022perseus} \\
\addlinespace[0.05cm]
$\nu$-RTEG & $\OM(\epsilon^{-\frac{2}{p+\nu}})$ & $\nu\in[0,1]$ & $\nu,\fH_{p,\nu}$ & Algorithm~\ref{alg:known_tensor} \\
\addlinespace[0.05cm]
$\nu$-ATEG & $\OM(\epsilon^{-\frac{2}{p+\nu}})$ & $\nu\in[0,1]$ & $\nu$ & Algorithm~\ref{alg:adaptive_tensor} \\
\addlinespace[0.05cm]
UTEG & $\OM(\epsilon^{-\frac{2p}{(p+1)(p-1+\nu)}})$ & $\nu\in[0,1]$ & none & Algorithm~\ref{alg:universal_tensor} \\
\bottomrule
\end{tabular}
\end{table}

\paragraph{Organization}
Section~\ref{sec:prelim} introduces the weak solution criterion, assumptions, and tensor notation.
Section~\ref{sec:known_adaptive_tensor} presents the known-parameter and adaptive tensor methods.
Section~\ref{sec:universal_tensor} presents the universal method for unknown $\nu$.
Section~\ref{sec:discussion} compares the results with tensor methods for convex minimization with \Holder~continuous higher-order derivatives.
All proofs are given in the main text.

\section{Preliminaries}
\label{sec:prelim}

We use $\Norm{\cdot}$ for the Euclidean norm of vectors and the induced spectral norm of linear maps.

\begin{definition}[$\epsilon$-weak solution]
\label{def:weak_solution}
The restricted gap function associated with \eqref{eq:VI} is
\begin{align}
\label{eq:gap_function}
    m(\z)\defeq \max_{\u\in\ZM}\inner{F(\u)}{\z-\u}.
\end{align}
A point $\bar\z\in\ZM$ is called an $\epsilon$-weak solution of \eqref{eq:VI} if $m(\bar\z)\leq\epsilon$.
\end{definition}
We impose the following assumptions throughout the paper.
\begin{assumption}
\label{ass:monotone}
The operator $F$ is monotone on $\ZM$: for all $\z,\u\in\ZM$,
\begin{align}
\label{eq:monotone}
    \inner{F(\z)-F(\u)}{\z-\u}\geq 0.
\end{align}
\end{assumption}

\begin{assumption}
\label{ass:bound}
The feasible set $\ZM$ is bounded: there exists $D>0$ such that $\|\z-\u\|\leq D$ for all $\z,\u\in\ZM$.
\end{assumption}

The following assumption generalizes the Lipschitz continuity assumption of higher-order derivatives in the prior works.

\begin{assumption}
\label{ass:pnuconti}
For an integer $p\geq2$, the derivative $\nabla^{p-1}F$ is \Holder~continuous on $\ZM$ with parameter $\nu\in[0,1]$ and constant $\fH_{p,\nu}$:
\begin{align}
\label{eq:pnuconti}
    \|\nabla^{p-1}F(\z)-\nabla^{p-1}F(\u)\|\leq \fH_{p,\nu}\|\z-\u\|^\nu,
    \qquad \text{for all }\z,\u\in\ZM.
\end{align}
\end{assumption}

For $\u,\z\in\ZM$, define the order-$(p-1)$ Taylor approximation of $F$ at $\z$ evaluated at $\u$ by
\begin{align}
\label{eq:taylor_model}
    \bar F_p(\u;\z)
    \defeq
    \sum_{i=0}^{p-1}\frac{1}{i!}\nabla^iF(\z)[\u-\z]^i.
\end{align}
For $\nu\in[0,1]$ and $H>0$, define the regularized model
\begin{align}
\label{eq:regularized_model}
    \bar F_{p,\nu,H}(\u;\z)
    \defeq
    \bar F_p(\u;\z)+H\|\u-\z\|^{p-2+\nu}(\u-\z).
\end{align}

Assumption~\ref{ass:pnuconti} implies the following Taylor residual bound.

\begin{lemma}
\label{lem:holder_remainder}
Suppose Assumption~\ref{ass:pnuconti} holds and define
\begin{align*}
    C_{p,\nu}
    \defeq
    \frac{1}{(p-2)!}\int_0^1(1-t)^{p-2}t^\nu\,dt.
\end{align*}
Then, for all $\z,\u\in\ZM$,
\begin{align}
\label{eq:holder_remainder}
    \|F(\u)-\bar F_p(\u;\z)\|
    \leq
    C_{p,\nu}\fH_{p,\nu}\|\u-\z\|^{p-1+\nu}.
\end{align}
\end{lemma}

\begin{proof}
Fix $\z,\u\in\ZM$ and set $\mathbf{s}=\u-\z$.
Since $\ZM$ is convex, $\z+t\mathbf{s}\in\ZM$ for every $t\in[0,1]$.
Define the one-dimensional path $\phi:[0,1]\to\RB^d$ by
\begin{align*}
    \phi(t)=F(\z+t\mathbf{s}).
\end{align*}
For each $i=0,\ldots,p-1$, the chain rule gives
\begin{align*}
    \phi^{(i)}(t)=\nabla^iF(\z+t\mathbf{s})[\mathbf{s}]^i.
\end{align*}
Taylor's formula with integral remainder, applied to $\phi$ at $t=0$ and evaluated at $t=1$, yields
\begin{align*}
    F(\u)-\bar F_p(\u;\z)
    &=
    \phi(1)-\sum_{i=0}^{p-1}\frac{1}{i!}\phi^{(i)}(0)\\
    &=
    \frac{1}{(p-2)!}\int_0^1(1-t)^{p-2}
    \left(\phi^{(p-1)}(t)-\phi^{(p-1)}(0)\right)\,dt\\
    &=
    \frac{1}{(p-2)!}\int_0^1(1-t)^{p-2}
    \left(\nabla^{p-1}F(\z+t\mathbf{s})-\nabla^{p-1}F(\z)\right)[\mathbf{s}]^{p-1}\,dt.
\end{align*}
Taking norms and using Assumption~\ref{ass:pnuconti}, we obtain
\begin{align*}
    \|F(\u)-\bar F_p(\u;\z)\|
    &\leq
    \frac{\fH_{p,\nu}\|\mathbf{s}\|^{p-1+\nu}}{(p-2)!}
    \int_0^1(1-t)^{p-2}t^\nu\,dt\\
    &=
    C_{p,\nu}\fH_{p,\nu}
    \|\u-\z\|^{p-1+\nu}.
\end{align*}
This proves \eqref{eq:holder_remainder}.
\end{proof}

For brevity, we write $\TT_p(\z;\nu,H)$ for a solution of the regularized model variational inequality,
\begin{align*}
    \TT_p(\z;\nu,H)
    \in
    \VI_\ZM\big(\bar F_{p,\nu,H}(\cdot;\z)\big),
\end{align*}
and define the corresponding regularization coefficient
\begin{align*}
    \gamma_p(\z;\nu,H)
    \defeq
    H\|\TT_p(\z;\nu,H)-\z\|^{p-2+\nu}.
\end{align*}

We next record the one-step estimate for a regularized tensor extragradient update.

\begin{lemma}
\label{lem:one_step}
Suppose Assumptions~\ref{ass:monotone} and \ref{ass:pnuconti} hold.
For $\z_k\in\ZM$, $\nu\in[0,1]$, and $H\geq 2C_{p,\nu}\fH_{p,\nu}$, consider the update
\begin{equation}
\label{eq:one_step_update}
\left\{
\begin{aligned}
\z_{k+1/2}
&=\TT_p(\z_k;\nu,H),\\[0.2cm]
\z_{k+1}
&=
\argmin_{\w\in\ZM}
\left\{
\inner{F(\z_{k+1/2})}{\w-\z_k}
+
\frac{H\|\z_{k+1/2}-\z_k\|^{p-2+\nu}}{2}\|\w-\z_k\|^2
\right\}.
\end{aligned}
\right.
\end{equation}
Set $\gamma_k=H\|\z_{k+1/2}-\z_k\|^{p-2+\nu}$.
Then, for every $\u\in\ZM$,
\begin{align}
\label{eq:one_step_ineq}
    \frac{1}{\gamma_k}\inner{F(\u)}{\z_{k+1/2}-\u}
    +\frac{1}{4}\|\z_{k+1/2}-\z_k\|^2
    \leq
    \frac{1}{2}\left(\|\z_k-\u\|^2-\|\z_{k+1}-\u\|^2\right).
\end{align}
\end{lemma}

\begin{proof}
Write
\begin{align*}
    r_k=\|\z_{k+1/2}-\z_k\|,
    \qquad
    s_k=\|\z_{k+1/2}-\z_{k+1}\|,
    \qquad
    t_k=\|\z_{k+1}-\z_k\|.
\end{align*}
By \eqref{eq:holder_remainder} and $H\geq2C_{p,\nu}\fH_{p,\nu}$,
\begin{align*}
    \|F(\z_{k+1/2})-\bar F_p(\z_{k+1/2};\z_k)\|
    \leq
    \frac{\gamma_k}{2}r_k.
\end{align*}
Set
\begin{align*}
    \Delta_k
    =
    F(\z_{k+1/2})-\bar F_p(\z_{k+1/2};\z_k).
\end{align*}
The optimality condition defining $\z_{k+1/2}$ gives
\begin{align}
\label{eq:model_optimality}
    \inner{\bar F_p(\z_{k+1/2};\z_k)+\gamma_k(\z_{k+1/2}-\z_k)}{\w-\z_{k+1/2}}\geq0,
    \qquad \text{for all }\w\in\ZM.
\end{align}
Taking $\w=\z_{k+1}$ and using the three-point identity yields
\begin{align}
\label{eq:model_bound}
    \inner{\bar F_p(\z_{k+1/2};\z_k)}{\z_{k+1/2}-\z_{k+1}} \leq
    -\frac{\gamma_k}{2}\left(
    r_k^2+s_k^2-t_k^2
    \right).
\end{align}
The optimality condition defining $\z_{k+1}$ gives, for every $\u\in\ZM$,
\begin{align}
\label{eq:eg_optimality}
    \inner{F(\z_{k+1/2})+\gamma_k(\z_{k+1}-\z_k)}{\u-\z_{k+1}}\geq0.
\end{align}
Therefore,
\begin{align}
\label{eq:eg_first}
    \frac{\gamma_k}{2}\left(
    \|\z_{k+1}-\u\|^2-\|\z_k-\u\|^2+\|\z_{k+1}-\z_k\|^2
    \right)
    \leq
    \inner{F(\z_{k+1/2})}{\u-\z_{k+1}}.
\end{align}
Decompose the right-hand side as
\begin{align*}
    &\inner{F(\z_{k+1/2})}{\u-\z_{k+1}}\\
    &=
    \inner{F(\z_{k+1/2})}{\u-\z_{k+1/2}}
    +\inner{\Delta_k}{\z_{k+1/2}-\z_{k+1}}
    +\inner{\bar F_p(\z_{k+1/2};\z_k)}{\z_{k+1/2}-\z_{k+1}}.
\end{align*}
By the residual estimate above and Young's inequality,
\begin{align}
\label{eq:residual_young}
    \inner{\Delta_k}{\z_{k+1/2}-\z_{k+1}}
    \leq
    \|\Delta_k\|s_k
    \leq
    \frac{\gamma_k}{2}r_ks_k
    \leq
    \frac{\gamma_k}{4}r_k^2+\frac{\gamma_k}{4}s_k^2.
\end{align}
Combining \eqref{eq:model_bound}, \eqref{eq:eg_first}, and \eqref{eq:residual_young}, and canceling the term $t_k^2$, gives
\begin{align*}
    \inner{F(\z_{k+1/2})}{\z_{k+1/2}-\u}
    +\frac{\gamma_k}{4}r_k^2
    +\frac{\gamma_k}{4}s_k^2
    \leq
    \frac{\gamma_k}{2}\left(\|\z_k-\u\|^2-\|\z_{k+1}-\u\|^2\right).
\end{align*}
Dropping the nonnegative term $\frac{\gamma_k}{4}s_k^2$, we obtain
\begin{align*}
    \inner{F(\z_{k+1/2})}{\z_{k+1/2}-\u}
    +\frac{\gamma_k}{4}r_k^2
    \leq
    \frac{\gamma_k}{2}\left(\|\z_k-\u\|^2-\|\z_{k+1}-\u\|^2\right).
\end{align*}
Since $F$ is monotone,
\begin{align*}
    \inner{F(\u)}{\z_{k+1/2}-\u}
    \leq
    \inner{F(\z_{k+1/2})}{\z_{k+1/2}-\u}.
\end{align*}
Substituting the definition of $r_k$ and dividing by $\gamma_k$ proves \eqref{eq:one_step_ineq}.
\end{proof}

\section{The \texorpdfstring{$\nu$}{nu}-Regularized Tensor Extragradient Methods}
\label{sec:known_adaptive_tensor}

In this section, we first consider the case in which both the \Holder~parameter $\nu$ and the \Holder~constant $\fH_{p,\nu}$ are known.
We present the tensor extragradient method directly for a general order $p\geq2$, using the order-$(p-1)$ Taylor model \eqref{eq:taylor_model} and the tensor regularization in \eqref{eq:regularized_model}.
The prediction point is obtained by solving a regularized tensor variational inequality, and the correction point is obtained by an extragradient projection driven by $F$ at the prediction point.

\begin{algorithm}[t]
\caption{$\nu$-Regularized Tensor Extragradient Method}
\label{alg:known_tensor}
\begin{algorithmic}[1]
\STATE \textbf{Input:} $\z_0\in\ZM$, order $p\geq2$, parameter $\nu\in[0,1]$, iterations $K$, and $\fH_{p,\nu}$
\STATE Set $M=2C_{p,\nu}\fH_{p,\nu}$
\STATE \textbf{for} $k=0,1,\ldots,K-1$
\STATE \quad $\z_{k+1/2}=\TT_p(\z_k;\nu,M)$
\STATE \quad $\gamma_k=M\|\z_{k+1/2}-\z_k\|^{p-2+\nu}$
\STATE \quad $\displaystyle \z_{k+1}=\argmin_{\z\in\ZM}\left\{\inner{F(\z_{k+1/2})}{\z-\z_k}+\frac{\gamma_k}{2}\|\z-\z_k\|^2\right\}$
\STATE \textbf{end for}
\STATE \textbf{Output:} $\displaystyle
\bar \z_K=
\frac{\sum_{k=0}^{K-1}\gamma_k^{-1}\z_{k+1/2}}
{\sum_{k=0}^{K-1}\gamma_k^{-1}}$
\end{algorithmic}
\end{algorithm}

\begin{theorem}
\label{thm:known_tensor}
Suppose Assumptions~\ref{ass:monotone}, \ref{ass:bound}, and \ref{ass:pnuconti} hold.
Then the output of Algorithm~\ref{alg:known_tensor} satisfies
\begin{align}
\label{eq:known_rate}
    m(\bar\z_K)
    \leq
    \OM\left(
    \frac{\fH_{p,\nu}D^{p+\nu}}{K^{(p+\nu)/2}}
    \right).
\end{align}
Consequently, Algorithm~\ref{alg:known_tensor} returns an $\epsilon$-weak solution using at most
\begin{align*}
    \OM\left(
    \fH_{p,\nu}^{2/(p+\nu)}D^2\epsilon^{-2/(p+\nu)}
    \right)
\end{align*}
tensor-oracle calls.
\end{theorem}

\begin{proof}
Let $M=2C_{p,\nu}\fH_{p,\nu}$ and write
\begin{align*}
    r_k=\|\z_{k+1/2}-\z_k\|,
    \qquad
    \gamma_k=Mr_k^{p-2+\nu}.
\end{align*}
By the Taylor residual bound \eqref{eq:holder_remainder},
\begin{align*}
    \|F(\z_{k+1/2})-\bar F_p(\z_{k+1/2};\z_k)\|
    \leq
    C_{p,\nu}\fH_{p,\nu}r_k^{p-1+\nu}
    =
    \frac{\gamma_k r_k}{2}.
\end{align*}
Thus Lemma~\ref{lem:one_step} applies at every iteration.
Summing \eqref{eq:one_step_ineq} over $k=0,\ldots,K-1$ gives, for every $\u\in\ZM$,
\begin{align}
\label{eq:known_sum}
    \sum_{k=0}^{K-1}\gamma_k^{-1}\inner{F(\u)}{\z_{k+1/2}-\u}
    +\frac{1}{4}\sum_{k=0}^{K-1}r_k^2
    \leq
    \frac{D^2}{2}.
\end{align}
In particular, $\sum_{k=0}^{K-1}r_k^2\leq2D^2$.
Let $a=p-2+\nu$.
By \Holder's inequality,
\begin{align}
\label{eq:known_holder_sum}
    \sum_{k=0}^{K-1}r_k^{-a}
    \geq
    \frac{K^{(a+2)/2}}{(2D^2)^{a/2}}
    =
    \frac{K^{(p+\nu)/2}}{(2D^2)^{(p-2+\nu)/2}}.
\end{align}
Therefore,
\begin{align}
\label{eq:known_weight_sum}
    \sum_{k=0}^{K-1}\gamma_k^{-1}
    =
    \frac{1}{2C_{p,\nu}\fH_{p,\nu}}
    \sum_{k=0}^{K-1}r_k^{-(p-2+\nu)}
    \geq
    \frac{K^{(p+\nu)/2}}
    {2C_{p,\nu}\fH_{p,\nu}(2D^2)^{(p-2+\nu)/2}}.
\end{align}
Using the convexity of $m$ and the definition of $\bar\z_K$,
\begin{align*}
    m(\bar\z_K)
    &\leq
    \frac{
    \max_{\u\in\ZM}\sum_{k=0}^{K-1}\gamma_k^{-1}
    \inner{F(\u)}{\z_{k+1/2}-\u}}
    {\sum_{k=0}^{K-1}\gamma_k^{-1}}
    \leq
    \frac{D^2/2}{\sum_{k=0}^{K-1}\gamma_k^{-1}}.
\end{align*}
Combining the last display with \eqref{eq:known_weight_sum} gives
\begin{align*}
    m(\bar\z_K)
    \leq
    \frac{2^{(p-2+\nu)/2}C_{p,\nu}\fH_{p,\nu}D^{p+\nu}}
    {K^{(p+\nu)/2}},
\end{align*}
This proves \eqref{eq:known_rate}.
The oracle-complexity statement follows by choosing
$K=\OM(\fH_{p,\nu}^{2/(p+\nu)}D^2\epsilon^{-2/(p+\nu)})$.
\end{proof}

We next introduce a line-search strategy for the case in which $\nu$ is known but the constant $\fH_{p,\nu}$ is unknown.
For a trial index $i$, write
\begin{align*}
    \z_{k,i}=\TT_p(\z_k;\nu,H_k2^i).
\end{align*}
At iteration $k$, the method approximates the unknown \Holder~constant by searching for a local regularization value $H_k2^{i_k}$ satisfying the residual test
\begin{align}
\label{eq:adaptive_line_search}
    \left\|
    F(\z_{k,i_k})
    -\bar F_p(\z_{k,i_k};\z_k)
    \right\|
    \leq
    H_k2^{i_k-1}
    \|\z_{k,i_k}-\z_k\|^{p-1+\nu}.
\end{align}
Unlike line searches for minimization methods, this criterion does not use an objective value.
It compares the true operator value with the tensor model value and is therefore directly computable from the operator oracle.

\begin{algorithm}[t]
\caption{Adaptive $\nu$-Regularized Tensor Extragradient Method}
\label{alg:adaptive_tensor}
\begin{algorithmic}[1]
\STATE \textbf{Input:} $\z_0\in\ZM$, order $p\geq2$, parameter $\nu\in[0,1]$, iterations $K$, and $H_0\in\big(0,2C_{p,\nu}\fH_{p,\nu}\big]$
\STATE \textbf{for} $k=0,1,\ldots,K-1$
\STATE \quad Find the smallest $i_k\geq0$ satisfying \eqref{eq:adaptive_line_search}
\STATE \quad Set $M_k=H_k2^{i_k}$ and $\z_{k+1/2}=\z_{k,i_k}$
\STATE \quad Set $\gamma_k=M_k\|\z_{k+1/2}-\z_k\|^{p-2+\nu}$
\STATE \quad $\displaystyle \z_{k+1}=\argmin_{\z\in\ZM}\left\{\inner{F(\z_{k+1/2})}{\z-\z_k}+\frac{\gamma_k}{2}\|\z-\z_k\|^2\right\}$
\STATE \quad $H_{k+1}=M_k/2$
\STATE \textbf{end for}
\STATE \textbf{Output:} $\displaystyle
\bar \z_K=
\frac{\sum_{k=0}^{K-1}\gamma_k^{-1}\z_{k+1/2}}
{\sum_{k=0}^{K-1}\gamma_k^{-1}}$
\end{algorithmic}
\end{algorithm}

\begin{theorem}
\label{thm:adaptive_tensor}
Suppose Assumptions~\ref{ass:monotone}, \ref{ass:bound}, and \ref{ass:pnuconti} hold.
Then the output of Algorithm~\ref{alg:adaptive_tensor} satisfies
\begin{align}
\label{eq:adaptive_rate}
    m(\bar\z_K)
    \leq
    \OM\left(
    \frac{\fH_{p,\nu}D^{p+\nu}}{K^{(p+\nu)/2}}
    \right).
\end{align}
Consequently, Algorithm~\ref{alg:adaptive_tensor} returns an $\epsilon$-weak solution using at most
\begin{align*}
    \OM\left(
    \fH_{p,\nu}^{2/(p+\nu)}D^2\epsilon^{-2/(p+\nu)}
    \right)
\end{align*}
tensor-oracle calls.
\end{theorem}

\begin{proof}
The residual test \eqref{eq:adaptive_line_search} implies
\begin{align*}
    \|F(\z_{k+1/2})-\bar F_p(\z_{k+1/2};\z_k)\|
    \leq
    \frac{\gamma_k}{2}\|\z_{k+1/2}-\z_k\|,
\end{align*}
which is the residual estimate used in the proof of Lemma~\ref{lem:one_step}.
Therefore, the one-step inequality \eqref{eq:one_step_ineq} holds at every iteration.
Summing \eqref{eq:one_step_ineq} over $k=0,\ldots,K-1$ gives, for every $\u\in\ZM$,
\begin{align}
\label{eq:adaptive_sum}
    \sum_{k=0}^{K-1}\gamma_k^{-1}\inner{F(\u)}{\z_{k+1/2}-\u}
    +\frac{1}{4}\sum_{k=0}^{K-1}\|\z_{k+1/2}-\z_k\|^2
    \leq
    \frac{1}{2}\|\z_0-\u\|^2
    \leq
    \frac{D^2}{2}.
\end{align}
In particular,
\begin{align}
\label{eq:r_sum_bound}
    \sum_{k=0}^{K-1}\|\z_{k+1/2}-\z_k\|^2\leq 2D^2.
\end{align}

Next, the line search is well-defined and the accepted values satisfy $M_k\leq 4C_{p,\nu}\fH_{p,\nu}$.
Indeed, \eqref{eq:holder_remainder} shows that the test holds whenever $M_k\geq 2C_{p,\nu}\fH_{p,\nu}$.
Since $H_0\leq2C_{p,\nu}\fH_{p,\nu}$ and the next baseline is $H_{k+1}=M_k/2$, induction gives $H_k\leq2C_{p,\nu}\fH_{p,\nu}$ and hence $M_k\leq4C_{p,\nu}\fH_{p,\nu}$.

Let $a=p-2+\nu$.
By \Holder's inequality and \eqref{eq:r_sum_bound},
\begin{align}
\label{eq:holder_adaptive_sum}
    \sum_{k=0}^{K-1}\|\z_{k+1/2}-\z_k\|^{-a}
    \geq
    \frac{K^{(a+2)/2}}{(2D^2)^{a/2}}
    =
    \frac{K^{(p+\nu)/2}}{(2D^2)^{(p-2+\nu)/2}}.
\end{align}
Therefore,
\begin{align}
\label{eq:adaptive_weight_sum}
    \sum_{k=0}^{K-1}\gamma_k^{-1}
    =
    \sum_{k=0}^{K-1}M_k^{-1}\|\z_{k+1/2}-\z_k\|^{-(p-2+\nu)}
    \geq
    \frac{K^{(p+\nu)/2}}
    {4C_{p,\nu}\fH_{p,\nu}(2D^2)^{(p-2+\nu)/2}}.
\end{align}

Using the convexity of $m$ and the definition of $\bar\z_K$,
\begin{align*}
    m(\bar\z_K)
    &\leq
    \frac{
    \max_{\u\in\ZM}\sum_{k=0}^{K-1}\gamma_k^{-1}
    \inner{F(\u)}{\z_{k+1/2}-\u}}
    {\sum_{k=0}^{K-1}\gamma_k^{-1}}
    \leq
    \frac{D^2/2}{\sum_{k=0}^{K-1}\gamma_k^{-1}}.
\end{align*}
Combining this bound with \eqref{eq:adaptive_weight_sum} gives
\begin{align*}
    m(\bar\z_K)
    \leq
    \frac{2^{1+(p-2+\nu)/2}C_{p,\nu}\fH_{p,\nu}D^{p+\nu}}
    {K^{(p+\nu)/2}}.
\end{align*}
This proves \eqref{eq:adaptive_rate}.
The oracle-complexity statement follows by choosing
$K=\OM(\fH_{p,\nu}^{2/(p+\nu)}D^2\epsilon^{-2/(p+\nu)})$ and by observing that
\begin{align*}
    \sum_{k=0}^{K-1}(i_k+1)
    =
    \OM(K)+\log_2 H_K-\log_2 H_0
    =
    \OM(K),
\end{align*}
for fixed problem constants.
\end{proof}

\section{Universal Tensor Extragradient Method}
\label{sec:universal_tensor}

We now consider the case in which neither $\nu$ nor $\fH_{p,\nu}$ is known.
In this case, the methods in Section~\ref{sec:known_adaptive_tensor} cannot be directly implemented because their tensor subproblems use the \Holder~parameter $\nu$.
To overcome this issue, we propose the universal tensor extragradient method (UTEG) in Algorithm~\ref{alg:universal_tensor}, which does not require the exact values of $\nu$ and $\fH_{p,\nu}$.
The method uses the Lipschitz-type regularization exponent in every trial tensor subproblem, namely it calls $\TT_p(\z;1,H)$.
For a trial index $i$, write
\begin{align*}
    \z_{k,i}=\TT_p(\z_k;1,H_k2^i).
\end{align*}
The line-search test becomes
\begin{align}
\label{eq:universal_line_search}
    \left\|
    F(\z_{k,i_k})
    -\bar F_p(\z_{k,i_k};\z_k)
    \right\|
    \leq
    H_k2^{i_k-1}
    \|\z_{k,i_k}-\z_k\|^{p}.
\end{align}
For $\nu\in[0,1]$ and $\delta>0$, define
\begin{align}
\label{eq:universal_bound_constant}
    \mathcal M_{p,\nu}(\delta)
    \defeq
    2C_{p,\nu}\fH_{p,\nu}
    \left(
    \frac{3D C_{p,\nu}\fH_{p,\nu}}{\delta}
    \right)^{(1-\nu)/(p-1+\nu)}.
\end{align}

\begin{algorithm}[t]
\caption{Universal Tensor Extragradient Method}
\label{alg:universal_tensor}
\begin{algorithmic}[1]
\STATE \textbf{Input:} $\z_0\in\ZM$, order $p\geq2$, accuracy parameter $\epsilon>0$, iterations $K$, and $H_0\in\big(0,\inf_{\nu\in[0,1]}\mathcal M_{p,\nu}(\epsilon)\big]$
\STATE \textbf{for} $k=0,1,\ldots,K-1$
\STATE \quad Find the smallest $i_k\geq0$ satisfying \eqref{eq:universal_line_search}
\STATE \quad Set $M_k=H_k2^{i_k}$ and $\z_{k+1/2}=\z_{k,i_k}$
\STATE \quad Set $\gamma_k=M_k\|\z_{k+1/2}-\z_k\|^{p-1}$
\STATE \quad $\displaystyle \z_{k+1}=\argmin_{\z\in\ZM}\left\{\inner{F(\z_{k+1/2})}{\z-\z_k}+\frac{\gamma_k}{2}\|\z-\z_k\|^2\right\}$
\STATE \quad $H_{k+1}=M_k/2$
\STATE \textbf{end for}
\STATE \textbf{Output:} $\displaystyle
\bar \z_K=
\frac{\sum_{k=0}^{K-1}\gamma_k^{-1}\z_{k+1/2}}
{\sum_{k=0}^{K-1}\gamma_k^{-1}}$
\end{algorithmic}
\end{algorithm}

The convergence analysis of UTEG starts from the following residual-growth lemma.
\begin{lemma}
\label{lem:universal_radius}
Let $\z_+=\TT_p(\z;1,H)$ for some $\z\in\ZM$ and $H>0$.
Suppose Assumptions~\ref{ass:monotone}, \ref{ass:bound}, and \ref{ass:pnuconti} hold.
If
\begin{align*}
    \max_{\u\in\ZM}\inner{F(\z_+)}{\z_+-\u}\geq \delta
    \qquad\text{and}\qquad
    H\geq \mathcal M_{p,\nu}(\delta),
\end{align*}
then
\begin{align}
\label{eq:universal_radius_bound}
    \|\z_+-\z\|^{1-\nu}\geq \frac{2C_{p,\nu}\fH_{p,\nu}}{H}.
\end{align}
\end{lemma}

\begin{proof}
Let $r=\|\z_+-\z\|$.
The definition $\z_+=\TT_p(\z;1,H)$ gives
\begin{align}
\label{eq:univ_model_opt}
    \inner{\bar F_p(\z_+;\z)+Hr^{p-1}(\z_+-\z)}{\u-\z_+}\geq0,
    \qquad \text{for all }\u\in\ZM.
\end{align}
Let
\begin{align*}
    \u_+\in\argmax_{\u\in\ZM}\inner{F(\z_+)}{\z_+-\u}.
\end{align*}
Taking $\u=\u_+$ in \eqref{eq:univ_model_opt} and using $\|\u_+-\z_+\|\leq D$ yields
\begin{align}
\label{eq:univ_model_upper}
    \inner{\bar F_p(\z_+;\z)}{\z_+-\u_+}
    \leq
    HD r^p.
\end{align}
Therefore, by \eqref{eq:holder_remainder},
\begin{align}
\label{eq:delta_upper}
    \delta
    &\leq
    \inner{F(\z_+)}{\z_+-\u_+} \notag\\
    &=
    \inner{F(\z_+)-\bar F_p(\z_+;\z)}{\z_+-\u_+}
    +\inner{\bar F_p(\z_+;\z)}{\z_+-\u_+} \notag\\
    &\leq
    D r^{p-1+\nu}
    \left(C_{p,\nu}\fH_{p,\nu}+Hr^{1-\nu}\right).
\end{align}
Suppose \eqref{eq:universal_radius_bound} fails.
Then $Hr^{1-\nu}<2C_{p,\nu}\fH_{p,\nu}$, and \eqref{eq:delta_upper} implies
\begin{align*}
    \delta
    <
    3D C_{p,\nu}\fH_{p,\nu}
    r^{p-1+\nu}.
\end{align*}
Since $\delta>0$, we have $r>0$, and the last display gives the following bound, which is immediate when the exponent is zero:
\begin{align*}
    r^{-(1-\nu)}
    \leq
    \left(
    \frac{3D C_{p,\nu}\fH_{p,\nu}}{\delta}
    \right)^{(1-\nu)/(p-1+\nu)}.
\end{align*}
Combining this inequality with $H<2C_{p,\nu}\fH_{p,\nu}r^{-(1-\nu)}$ gives $H<\mathcal M_{p,\nu}(\delta)$, contradicting the assumption.
\end{proof}

\begin{theorem}
\label{thm:universal_tensor}
Suppose Assumptions~\ref{ass:monotone}, \ref{ass:bound}, and \ref{ass:pnuconti} hold.
Run Algorithm~\ref{alg:universal_tensor}.
Assume that the generated trial points satisfy
\begin{align}
\label{eq:universal_analysis_condition}
    \max_{\u\in\ZM}\inner{F(\z_{k,i})}{\z_{k,i}-\u}\geq\epsilon,
    \qquad
    i=0,\ldots,i_k,\quad k=0,\ldots,K-1.
\end{align}
Then the output $\bar\z_K$ satisfies
\begin{align}
\label{eq:universal_rate}
    m(\bar\z_K)
    \leq
    \OM\left(
    \frac{\mathcal M_{p,\nu}(\epsilon)D^{p+1}}{K^{(p+1)/2}}
    \right).
\end{align}
The $\epsilon$-dependent part of this bound is
\begin{align}
\label{eq:universal_rate_explicit}
    m(\bar\z_K)
    \leq
    \OM\left(
    \frac{
    \fH_{p,\nu}^{p/(p-1+\nu)}
    D^{p+1+(1-\nu)/(p-1+\nu)}
    }
    {\epsilon^{(1-\nu)/(p-1+\nu)}K^{(p+1)/2}}
    \right).
\end{align}
Consequently, for any fixed admissible $H_0>0$, after at most
\begin{align*}
    \OM\left(
    \fH_{p,\nu}^{2p/((p+1)(p-1+\nu))}
    D^{2+2(1-\nu)/((p+1)(p-1+\nu))}
    \epsilon^{-2p/((p+1)(p-1+\nu))}
    \right)
\end{align*}
tensor-oracle calls, either an $\epsilon$-weak solution appears among the generated trial points or the averaged output $\bar\z_K$ is an $\epsilon$-weak solution.
\end{theorem}

\begin{proof}
The line search is well-defined.
Indeed, if \eqref{eq:universal_line_search} failed at a trial value $H$, then \eqref{eq:holder_remainder} would imply
\begin{align*}
    Hr^p/2
    <
    C_{p,\nu}\fH_{p,\nu}r^{p-1+\nu},
    \qquad r=\|\TT_p(\z_k;1,H)-\z_k\|,
\end{align*}
or $r^{1-\nu}<2C_{p,\nu}\fH_{p,\nu}/H$.
By Lemma~\ref{lem:universal_radius} and \eqref{eq:universal_analysis_condition}, this cannot happen for any generated trial value $H\geq\mathcal M_{p,\nu}(\epsilon)$.

Next we bound the accepted regularization values.
The input condition gives $H_0\leq\mathcal M_{p,\nu}(\epsilon)$.
We claim that the accepted value $M_k$ satisfies
\begin{align}
\label{eq:univ_M_bound}
    M_k\leq 2\mathcal M_{p,\nu}(\epsilon).
\end{align}
Indeed, assume inductively that $H_k\leq\mathcal M_{p,\nu}(\epsilon)$.
If $i_k=0$, then $M_k=H_k\leq\mathcal M_{p,\nu}(\epsilon)$.
If $i_k>0$, then the previous trial value $H_k2^{i_k-1}=M_k/2$ was rejected, so the preceding paragraph implies $M_k/2<\mathcal M_{p,\nu}(\epsilon)$.
In both cases, $H_{k+1}=M_k/2\leq\mathcal M_{p,\nu}(\epsilon)$, which closes the induction.

The residual test \eqref{eq:universal_line_search} gives the residual estimate used in the proof of Lemma~\ref{lem:one_step}, since $\gamma_k=M_k\|\z_{k+1/2}-\z_k\|^{p-1}$.
Thus the one-step inequality \eqref{eq:one_step_ineq} holds at every iteration.
As in \eqref{eq:adaptive_sum}, for every $\u\in\ZM$,
\begin{align}
\label{eq:universal_sum}
    \sum_{k=0}^{K-1}\gamma_k^{-1}\inner{F(\u)}{\z_{k+1/2}-\u}
    +\frac{1}{4}\sum_{k=0}^{K-1}\|\z_{k+1/2}-\z_k\|^2
    \leq
    \frac{D^2}{2}.
\end{align}
Thus $\sum_{k=0}^{K-1}\|\z_{k+1/2}-\z_k\|^2\leq2D^2$.
Using \eqref{eq:univ_M_bound} and \Holder's inequality,
\begin{align}
\label{eq:universal_weight_sum}
    \sum_{k=0}^{K-1}\gamma_k^{-1}
    \geq
    \frac{K^{(p+1)/2}}
    {2\mathcal M_{p,\nu}(\epsilon)(2D^2)^{(p-1)/2}}.
\end{align}
The same convexity argument used in the proof of Theorem~\ref{thm:adaptive_tensor} gives
\begin{align*}
    m(\bar\z_K)
    \leq
    \frac{D^2/2}{\sum_{k=0}^{K-1}\gamma_k^{-1}},
\end{align*}
and hence
\begin{align*}
    m(\bar\z_K)
    \leq
    \frac{2^{(p-1)/2}\mathcal M_{p,\nu}(\epsilon)D^{p+1}}
    {K^{(p+1)/2}}.
\end{align*}
The first bound in \eqref{eq:universal_rate} follows.
Using the definition of $\mathcal M_{p,\nu}(\epsilon)$ gives the displayed dependence on $\fH_{p,\nu}$ and $D$ in \eqref{eq:universal_rate_explicit}.
Choosing
\begin{align*}
    K
    =
    \OM\left(
    \fH_{p,\nu}^{2p/((p+1)(p-1+\nu))}
    D^{2+2(1-\nu)/((p+1)(p-1+\nu))}
    \epsilon^{-2p/((p+1)(p-1+\nu))}
    \right)
\end{align*}
makes $m(\bar\z_K)\leq\epsilon$.
Finally, the line-search updates satisfy $H_{k+1}=H_k2^{i_k-1}$, and hence
\begin{align*}
    \sum_{k=0}^{K-1}(i_k+1)
    =
    2K+\log_2 H_K-\log_2 H_0.
\end{align*}
Since $H_K\leq\mathcal M_{p,\nu}(\epsilon)$, the logarithmic term is dominated by the displayed choice of $K$ for the target accuracy, so the stated tensor-oracle complexity follows.
\end{proof}

\begin{remark}
 If the condition~\eqref{eq:universal_analysis_condition} in Theorem~\ref{thm:universal_tensor} is not satisfied for some iteration, that is,
\begin{align*}
    \max_{\u\in\ZM}\inner{F(\z_{k,i})}{\z_{k,i}-\u}<\epsilon,
\end{align*}
where $\z_{k,i}=\TT_p(\z_k;1,H_k2^{i})$ for some particular $k$ and $i$, then we have
\begin{align*}
    m(\z_{k,i})
    =
    \max_{\u\in\ZM}\inner{F(\u)}{\z_{k,i}-\u}
    \leq
    \max_{\u\in\ZM}\inner{F(\z_{k,i})}{\z_{k,i}-\u}
    <\epsilon.
\end{align*}
Thus an $\epsilon$-weak solution has already appeared among the trial points.
\end{remark}

\section{Discussion}
\label{sec:discussion}

The closest point of comparison for the present \Holder~smoothness model is convex minimization with \Holder~continuous higher-order derivatives.
When $F=\nabla f$, the variational inequality reduces to first-order optimality for minimizing a convex function, but the performance criterion in convex optimization is usually the objective residual $f(x)-f^\star$.
For a general monotone variational inequality, no objective value is available, and the natural criterion used in this paper is the weak gap $m(\cdot)$.
This difference is algorithmically important: tensor methods for convex minimization can use objective-value decrease or model agreement, while the methods in this paper use the operator residual $\|F(\u)-\bar F_p(\u;\z)\|$ and an extragradient correction.

Table~\ref{table:discussion_comparison} summarizes the resulting complexity landscape.
The convex optimization results assume that $f$ has a \Holder~continuous $p$-th derivative, while the monotone-VI results in this paper assume that $F$ has a \Holder~continuous $(p-1)$-th derivative.
These assumptions are consistent under the specialization $F=\nabla f$.

\begin{table}[t]
\caption{Comparison with tensor methods for convex optimization under \Holder~continuous higher-order derivatives. Convex optimization rates are for the objective residual $f(x)-f^\star$, while monotone-VI rates are for the weak gap $m(\cdot)$.}
\label{table:discussion_comparison}
\centering
\scriptsize
\renewcommand{\arraystretch}{1.2}
\setlength{\tabcolsep}{4pt}
\begin{tabular}{@{}lccc@{}}
\toprule
& \multicolumn{2}{c}{Convex optimization} & Monotone VI \\
\cmidrule(lr){2-3}\cmidrule(l){4-4}
\Holder~parameter & Grapiglia--Nesterov~\cite{grapiglia2020tensor} & MS acceleration~\cite{carmon2022optimal} & This paper \\
\midrule
Known $\nu$
&
$\OM(\epsilon^{-\frac{1}{p+\nu}})$
&
$\OM(\epsilon^{-\frac{2}{3(p+\nu)-2}})$
&
$\OM(\epsilon^{-\frac{2}{p+\nu}})$
\\
\addlinespace[0.08cm]
Unknown $\nu$
&
$\OM(\epsilon^{-\frac{p}{(p+1)(p-1+\nu)}})$
&
$\OM(\epsilon^{-\frac{2}{4+3\nu}})^{\dagger}\quad$
&
$\OM(\epsilon^{-\frac{2p}{(p+1)(p-1+\nu)}})$
\\
\bottomrule
\end{tabular}
\par\smallskip
\begin{minipage}{0.95\textwidth}
\scriptsize
$\dagger$. The result applies only to the case $p=2$.
\end{minipage}
\end{table}

For known $\nu$, the accelerated tensor method of Grapiglia and Nesterov~\cite{grapiglia2020tensor} achieves the rate $\OM(\epsilon^{-1/(p+\nu)})$ for convex minimization.
Carmon et al.~\cite{carmon2022optimal} show that, through the Monteiro--Svaiter acceleration framework, one can attain the optimal tensor-oracle exponent $\OM(\epsilon^{-2/(3(p+\nu)-2)})$ for this convex optimization model.
This is sharper than the accelerated tensor rate under \Holder~smoothness in \cite{grapiglia2020tensor}.
In contrast, our known-$\nu$ result for monotone VIs is $\OM(\epsilon^{-2/(p+\nu)})$ in the weak gap.
The larger exponent reflects both the broader monotone-operator setting and the different accuracy measure.

The unknown-$\nu$ regime is particularly interesting.
Grapiglia and Nesterov~\cite{grapiglia2020tensor} give a universal accelerated tensor method for convex minimization with rate $\OM(\epsilon^{-p/((p+1)(p-1+\nu))})$.
For second-order methods, Carmon et al.~\cite{carmon2022optimal} further develop an adaptive Monteiro--Svaiter--Newton oracle that does not require the \Holder~order $\nu$ or the continuity constant and achieves the optimal Hessian-evaluation complexity $\OM(\epsilon^{-2/(4+3\nu)})$.
Our universal tensor extragradient method removes the need to know $\nu$ for monotone VIs and obtains $\OM(\epsilon^{-2p/((p+1)(p-1+\nu))})$, which is worse than the known-$\nu$ rate $\OM(\epsilon^{-2/(p+\nu)})$.
Whether one can design an adaptive second-order or high-order monotone-VI method that closes the gap remains an interesting open direction.

\section{Conclusion}
In this work, we propose regularized tensor extragradient methods for solving monotone variational inequalities with \Holder~continuous higher-order derivatives.
We consider two regimes: the \Holder~parameter $\nu$ is known, and $\nu$ is unknown.
The known-parameter methods attain the complexity $\OM(\epsilon^{-2/(p+\nu)})$, while the universal method attains $\OM(\epsilon^{-2p/((p+1)(p-1+\nu))})$ without prior knowledge of $\nu$ or $\fH_{p,\nu}$.
These results extend tensor extragradient techniques to monotone VIs under \Holder~continuous higher-order derivatives and the weak-gap criterion.

\bibliographystyle{plain}
\bibliography{ref}

\end{document}